\documentclass[a4paper,12pt,reqno]{amsart}
\usepackage{amsmath,amsfonts,amssymb,amsthm}%

\numberwithin{equation}{section}
\newtheorem{thm}{Theorem}[section]

\newtheorem{lemma}[thm]{Lemma}

\newtheorem{cor}[thm]{Corollary}

\newtheorem{prop}[thm]{Proposition}

\newtheorem{rem}[thm]{Remark}

\newtheorem*{ack}{Acknowledgement}

\begin{document}
%\begin{titlepage}
\title{Analytic Dirac approximation for real 
linear algebraic groups}
\author{Christoph Lienau} 
%\thanks{Former adress: Max-Planck-Institut f\"ur Mathematik, Vivatsgasse 7, D-53111 Bonn, lienau@mpim-bonn.mpg.de}
%\maketitle
%\thispagestyle{empty}
%\noindent
%\address{\textsc{Leibniz Universit\"at Hannover}\\
%\textsc{Institut f\"ur Analysis}\\  
%\textsc{Welfengarten 1} \\ \textsc{D-30167 Bonn}
%\\\textsc{email: lienau@math.uni-hannover.de} \\
%\textsc{tel: 0049511/7623514}}
%\end{titlepage}
\title[Analytic Dirac approximation]{Analytic Dirac approximation for real linear algebraic groups}
\maketitle
\begin{abstract}
For a real linear algebraic group $G$ let $\mathcal{A}(G)$ be the algebra
of analytic vectors for the left regular representation of $G$ on the space of 
superexponentially decreasing functions. We present an explicit Dirac sequence
in $\mathcal{A}(G)$. Since $\mathcal{A}(G)$ acts on $E$ for every Fr\'{e}chet-representation $(\pi,E)$ of moderate growth, this yields an elementary proof of a result of Nelson that the space of analytic vectors is dense in $E$.
\end{abstract}\par
\small{\textsc{Subject Classification:} 22E30
\section{introduction}
\noindent
In this paper we provide an explicit Dirac sequence of superexponentially decrasing analytic
functions on a linear algebraic group. This yields an elementary proof of a theorem of Nelson \cite{N}
that the space of analytic vectors is dense. In order to keep the exposition self contained we 
recall basic constructions from \cite{GKL,GKS}\\
Let $(\pi,E)$ be a representation of a Lie group $G$ on a Fr\'{e}chet-space E. 
For a vector $v\in E$ we denote by $\gamma_{v}$ the corresponding orbit map
\[\gamma_{v}: \ G\rightarrow E, \ \ \ \  g\mapsto \pi(g)v.\]
A vector $v\in E$ is called {\it analytic} if the orbit map $\gamma_{v}$ is a real analytic $E$-valued map. 
We denote the space of all analytic vectors by $E^{\omega}$.\\
Let $\bold{g}$ be a left invariant Riemannian measure on $G$.
To $\bold{g}$ we associate a Riemannian distance $d$ on $G$: $d(g)$ is defined as the infimum lenght
of all arcs joining $g$ and $1$. \\
Let $\mathcal{R}(G)$ be the space of superexponentially decrasing smooth functions on $G$ with 
respect to the distance $d$, i.e 
\[\mathcal{R}(G)=\left\{f\in C(G) \ | \ \forall n\in \mathbb{N}: p_{n}(f):=\sup_{g\in G}|f(g)|e^{nd(g)}<\infty \right\}.
\]
The space $\mathcal{R}(G)$ is a Fr\'{e}chet algebra under convolution and is independent of the choice
of the left invariant metric. \\
Let us assume that $(\pi,E)$ is a $F$-representation, i.e the representation is of moderate growth. In particular every
Banach representation is a $F$-representation. \\
For every continuous seminorm $q$ on $E$ exists a continuous seminorm $q^{\prime}$ and constants $C,c>0$ such that
\[q\left(\pi(g)v\right)\leq Ce^{cd(g)}q^{\prime}(v) \ \ \ (\forall g\in G,\forall v\in E).
\]
Furthermore there exists a constant $r^{\prime}>0$ such that $\forall r>r^{\prime}$
\[\int_{G}e^{-rd(g)}\ dg<\infty.
\]
Hence there is a  corresponding algebra representation $\Pi$ of $\mathcal{R}(G)$ which is given by
\[\Pi(f)v=\int_{G}f(g)\pi(g)v\ dg \ \ \ (f\in\mathcal{}(G),v \in E).
\]
We denote the space of analytic vectors $\mathcal{R}(G)^{\omega}$ for the left regular representation  $L$ by $\mathcal{A}(G)$. \\
A function $f\in \mathcal{R}(G)$ is in $\mathcal{A}(G)$ if and only if it satisfies
the following two conditions.
\begin{enumerate}
\item \label{c1} There exists a neighborhood $U\subset G_{\mathbb{C}}$ of $1$ and a $F\in \mathcal{O}(U^{-1}G)$
with $F|_{G}=f$. 
\item \label{c2} For every compact subset $Q\subset U$ we have \\
$\sup_{k\in Q}p_{n}\left(L_{k}(F)\right)<\infty$ for all $n\in\mathbb{N}$.
\end{enumerate}
Throughout this text we refer to these conditions as condition $(\ref{c1})$ and condition 
$(\ref{c2})$.\\
We define a positive function on $GL_{n}(\mathbb{R})$ by 
\[
|g|=\sqrt{\operatorname{tr}\left(g^{t}g\right )} \ \ \  (g\in GL_{n}(\mathbb{R})).
\]
Let $K$ be the maximal compact subgroup $O(n)$ of $GL_{n}(\mathbb{R})$ and $K_{\mathbb{C}}$ its complexification. Then $|\cdot|$
is $K_{\mathbb{C}}$-bi-invariant and sub-multiplicative. Note that for a matrix 
$g=\left(a_{ij}\right)_{1\leq i,j\leq n}$ we have $|g|=\sqrt{\sum_{1\leq i,j\leq n}a_{ij}^{2}}$.
Hence $|\cdot|^{2}$ is holomorphic on $GL_{n}(\mathbb{C})$.\\
Let $G$ be a real linear algebraic group, then $G$ is a closed subgroup of some $GL_{n}(\mathbb{R})$.  \\
We define a norm in the sense of \cite{Wal} on $G$ by
\[\|g\|=\max\{|g|,|g^{-1}|\}, \ \ \ \ \left(g\in G\right).
\]
For $t>0$ we consider the function
\[\varphi_{t}: G  \rightarrow \mathbb{R}, \ \ g  \mapsto C_{t}e^{-t^{2}\left(|g-1|^{4n}+|g^{-1}-1|^{4n}\right)},
\]
with constants $C_{t}>0$ such that $\|\varphi_{t}\|_{\mathrm{L}^{1}(G)}=1$. \\
Recall that a sequence $\left(f_{k}\right)_{k>0}$ is called a Dirac
sequence  if it satisfies the following three conditions:

\renewcommand{\theenumi}{\alph{enumi}}
\begin{enumerate}
\item $f\geq 0$, $\forall k\in \mathbb{N}$
\item $\int_{G}f_{k}(g)\ dg=1$ , $\forall k\in \mathbb{N}$
\item For every $\varepsilon >0$ and every neighborhood $U$ of $1$ in $G$ exists a $M\geq1$ such that $\int_{G\backslash U}f_{m}(g)\ dg<\varepsilon$, $\forall m\geq M$.
\end{enumerate}
We prove the following theorem
\begin{thm}\label{t1}
\begin{enumerate}
\item $\varphi_{t}\in \mathcal{A}(G)$ for all $t>0$. 
\item The sequence $\left(\varphi_{t}\right)_{t>0}$ forms a Dirac sequence.
\end{enumerate}
\end{thm}
As a corollary we obtain a result of Nelson \cite{N} for real linear algebraic
groups.
\begin{cor}
Let $\left(\pi,E\right)$ be a $F$-representation of a real linear algebraic group $G$ on a
Fr\'{e}chet space $E$, then the space of analytic vectors $E^{\omega}$ is dense in $E$.
\end{cor}
\begin{rem}
In fact \cite{GKL} every analytic vector is a finite sum of vectors of the form
$\Pi(f)v$ with $f\in\mathcal{A}(G)$ and $v\in E$.  
\end{rem}
\begin{rem}
If $G$ is a real reductive group then Theorem \ref{t1} holds even for 
\[\varphi_{t}^{\prime}: G  \rightarrow \mathbb{R}, \ \ g  \mapsto C_{t}e^{-t^{2}\left(|g-1|^{2}+|g^{-1}-1|^{2}\right)},
\]
with constants $C_{t}>0$ such that $\|\varphi_{t}^{\prime}\|_{\mathrm{L}^{1}(G)}=1$.
\end{rem}
\section{Proofs}
The function $\varphi_{t}$ posses an holomorphic continuation to $G_{\mathbb{C}}$ which we also denote by
$\varphi_{t}$, but $\varphi_{t}$ does not satisfy condition $(2)$ on the whole of $G_{\mathbb{C}}$. \\
We now describe for $GL_{n}(\mathbb{R})$ a $K_{\mathbb{C}}\times GL_{n}(\mathbb{R})$-invariant
domain in $GL_{n}(\mathbb{C})$ where $\varphi_{t}$ satisfies condition $(2)$.
It turns out that this domain is a subdomain of the {\it crown domain} $\Xi$ \cite{AG, KS}. Therefore let 
$\Omega=\{\operatorname{diag}(d_{1},\ldots,d_{n})\  : \ d_{k}\in\mathbb{R},|d_{k}|<\frac{\pi}{4}, \forall \ k=1\ldots,n \}\subset\mathbb{R}^{n^{2}}$.
\begin{rem} Note that this $\Omega$ is not the same as in \cite{KS}. Let us denote by $\Omega_{ss}$ the
Omega used in \cite{KS} for $SL_{n}(\mathbb{R})$. Then $\Omega$ is related to $\Omega_{ss}$ in the following way: $\Omega$ has the property that $\Omega_{ss}+\mathbb{R}e=\Omega+\mathbb{R}e$ with
$e=\operatorname{diag}(1,\ldots,1)$. In other words, up to central shift the Omegas coincide.
\end{rem}
We define $\Xi_{n}$ by $\Xi_{n}=GL_{n}(\mathbb{R})\exp\left(i\frac{1}{n+1}\Omega\right)K_{\mathbb{C}}$.\\
\begin{rem}\label{rr}
Let us remark that if $G$ be a real reductive group, i.e a closed transposition stable subgroup
of $GL_{n}(\mathbb{R})$, then $d(g)$ and $\log\|g\|$ are comparable in the sense that there are constants
$c_{1},c_{2}>0$ and $C_{1},C_{2}\in \mathbb{R}$ such that
\[c_{1}d(g)+C_{1}\leq \log {\|g\|}\leq c_{1}d(g)+C_{2}.\]
Hence we can give an alternative characterization of the space $\mathcal{R}(G)$ in terms of $\|\cdot\|$:
\[\mathcal{R}(G)=\left\{f\in C(G) \ | \ \forall n\in \mathbb{N}: p_{n}(f):=\sup_{g\in G}\|g\|^{n}|f(g)|<\infty \right\}.
\]
\end{rem}
In the proof of the next proposition we need the following notations:
For a matrix $g=\left(a_{ij}\right)_{1\leq i,j\leq n}$ we denote by $g_{i}$ the i-th
column vector $(a_{1i},\ldots,a_{ni})^{T}$ and for a vector $w\in \mathbb{C}^{n}$ we denote by $\|w\|_{2}$ the euclidean norm. 
\begin{prop}\label{cd}
The function $\varphi_{t}$ satisfies condition $(\ref{c2})$ on $\Xi_{n}$.
\end{prop}
\begin{proof}
Let $Q\subset \Xi_{n}$ be compact. We show that there exists a constant $C>0$ such that 
\begin{align}\label{fs}
|\varphi_{t}(gq)|\leq e^{-C\|g\|^{4n}} , \ \ \ \ (\forall g\in G, \forall q\in Q).
\end{align}
There exists a $\Omega^{\prime}\subset \Omega$ which satisfies the 
following properties. 
\begin{enumerate}
\item $Q\subset GL_{n}(\mathbb{R})\exp(i\frac{1}{n+1}\Omega^{\prime})K_{\mathbb{C}}$. 
\item There exists a constant $C_{1}^{\prime}>0$ such that for all \\ $d=\operatorname{diag}(e^{i\theta_{1}},\ldots,e^{i\theta_{n}})\in \exp(i\frac{1}{n+1}\Omega^{\prime})$ we have \\ $\cos(2(\theta_{\alpha_{1}}+\ldots+\theta_{\alpha_{n+1}}))\geq C_{1}^{\prime}$ for all $\alpha_{j}\in \{1,\ldots,2n\}$. 
\end{enumerate}
This implies that for $k=1,\ldots,2n$ there exists a constant $C_{1}>0$ such that 
\begin{align}\label{p1}\operatorname{Re}\left(\left|gq\right|^{2k}\right)\geq C_{1} |g|^{2k}, \ \ \  \left(g\in GL_{n}(\mathbb{R}),q\in Q\right).
\end{align}
Therefore let $d=\operatorname{diag}(e^{i\theta_{1}},\ldots,e^{i\theta_{n}})\in \exp(i\Omega^{\prime})$ and $g^{\prime}\in GL_{n}(\mathbb{R})$ then
\begin{align}
|g^{\prime}d|^{2k}&=\left(e^{2\theta_{1}i}\|g_{1}^{\prime}\|^{2}+\cdots+e^{2\theta_{n}i}\|g_{n}^{\prime}\|^{2}\right)^{k}
\end{align}
Hence $\operatorname{Re}\left(\left|g^{\prime}q\right|^{2k}\right)\geq C_{1}^{\prime}|g^{\prime}|^{2k}$ according to $(b)$.\\
Let $q=hdk$ with $h\in GL_{n}(\mathbb{R})$ and $k\in K_{\mathbb{C}}$.  \\Then
$\operatorname{Re}\left(|gq|^{2}\right)=\operatorname{Re}\left(|ghdk|^{2}\right)=\operatorname{Re}\left(|ghd|^{2}\right)\geq C_{1}^{\prime}|gh|^{2}$. Since $Q$ is compact there exists a constant $C_{1}>0$ 
such that $C_{1}^{\prime}|gh|^{2}>C_{1}|g|^{2}$ for all $q\in Q$. Thus we obtain $(\ref{p1})$.\\
Likewise we can show that for $k=1,\ldots,2n$ there exists a constant $C_{2}>0$ such that
\begin{align}\label{p2}
\operatorname{Re}\left(\left|\left(gq\right)^{-1}\right|^{2k}\right)\geq C_{2} |g^{-1}|^{2k}, \ \ \  \left(g\in GL_{n}(\mathbb{R}),q\in Q\right).
\end{align}
Note that for $k=1,\ldots,4n$ there exists a constant $C_{3}>0$ such that 
\begin{align}\label{u2}
\operatorname{Re}\left(\operatorname{tr}\left(gq\right)^{k}\right)\leq |\operatorname{tr}\left(gq\right)|^{k} \leq 
 C_{3}|g|^{k} ,  \ \left(g\in GL_{n}(\mathbb{R}),q\in Q\right). 
\end{align}
Since $|g-1|^{4n}=(|g|^{2}-2\operatorname{tr}(g)+n)^{2n}$ we obtain the upper bound $(\ref{fs})$ for some $C>0$ by expanding the $2n$-th power and combining the estimates. \\
Since $GL_{n}(\mathbb{R})$ is real reductive Remark $\ref{rr}$ implies that $\varphi_{t}$ satisfies condition $(\ref{c2})$ on $\Xi$.
\end{proof}
Hence $\varphi_{t}\in \mathcal{A}\left(GL_{n}(\mathbb{R})\right)$. Now we show that
for every real linear algebraic group $G\subset GL_{n}(\mathbb{R})$ the functions $\varphi_{t}$ are elements
of $\mathcal{A}(G)$.
\begin{prop}
Let $G\subset GL_{n}(\mathbb{R})$ be a real linear algebraic group then $\varphi_{t}\in \mathcal{A}(G)$.
\end{prop}
\begin{proof}
The set $G_{\mathbb{C}}\cap \Xi_{n}$ is an open neighborhood of $1\in G$ to which $\varphi_{t}$ extends holomorphically.\\
We give an upper bound for $d(g)$ which implies that $\varphi_{t}$ satisfies $(\ref{c2})$ on this neighborhood.\\
Every algebraic group $G$ can be decomposed as a semidirect product $G=\operatorname{Rad}_{u}G\rtimes L$ of a connected unipotent group $\operatorname{Rad}_{u}G$
and a reductive group $L$. We write $g=ur$ with $u$ unipotent and $r$ reductive, then $d(g)=d(ur)=d(u)+d(r)$.\\
Remark \ref{rr} implies that there exists a constant $C>0$ such that $d(r)\leq C \log(\|r\|)+C$. 
Note that the unipotent radical $\operatorname{Rad}_{u}G$ is connected and 
$u$ has a real logarithm. The path $\gamma(t)=\exp(t\log(u))$ connects $1$ and $\log(u)$ and has length $|\log(u)|$, thus $d(u)\leq |\log(u)|$. Since $\log(u)=\sum_{k=0}^{n}\frac{(-1)^{k}(u-1)^{k}}{k}$ and $|u-1|^{k}\leq 1+|u-1|^{n}\leq 1+|u|^{n}$  for $k=0,\ldots,n$ we obtain $|\log(u)|\leq 1+n+n\|u\|^{n}$. 
Let $J=D+N$ be the Jordan normal form of $g$ with $D$ a diagonal and $N$ a nilpotent matrix and let $P\in GL_{n}(\mathbb{C})$ be the change of basis matrix. Since the Jordan-Chevalley
decomposition is unique, $u=P(1+D^{-1}N)P^{-1}$ and $r=PDP^{-1}$. Therefore $\|u\|\leq \|P\|^{2}(\|1\|+\|D^{-1}N\|)\leq\|P\|^{2}(\|1\|+\|D\|)\leq\|P\|^{2}(\|1\|+\|g\|)$. The last inequality follows from the
fact that the sum of the absolute values of the squares of the eigenvalues is less or equal than the sum of the squares of the singular values. Likewise we obtain $\|r\|\leq \|P\|^{2}\|g\|$. Since the column vectors of the matrices $P$ and $P^{-1}$ are chains of
generalized eigenvectors of $g$ we obtain $\|P\|^{2}\leq n2^{n}\|g\|^{2}$. \\
Combining these estimates we obtain that there exists a constant $R>0$ such that
\begin{align*}
e^{nd(g)}\leq Re^{R\|g\|^{3n}}\|g\|^{R}.
\end{align*}
Hence $\varphi_{t}$ satisfies condition $(\ref{c2})$ on $G_{\mathbb{C}}\cap\Xi_{n}$.
\end{proof}
\begin{prop}\label{d} The family $(\varphi_{t})_{t\geq1}$ forms for $t\rightarrow \infty$ a Dirac sequence. 
\end{prop}
\begin{proof}
Let $V$ be a neighborhood of $0$ in $\frak{g}$ such that the exponential map is a diffeomorphism
of $V$ with some neighborhood $U$ of $1$ in $G$. Then
\begin{align*}
\int_{G}e^{-t^{2}\left(|g-1|^{2}+|g^{-1}-1|^{2}\right)}\ dg\geq \int_{U}e^{-t^{2}\left(|g-1|^{2}+|g^{-1}-1|^{2}\right)}\ dg.
\end{align*}
The differential of $\operatorname{exp}$ at $X$ is given by
\[dL_{\exp(X)}\circ\tfrac{1-e^{-\operatorname{ad}(X)}}{\operatorname{ad}(X)}.
\]
Therefore
\begin{align*}
\int_{U}e^{-t^{2}\left(|g-1|^{2}+|g^{-1}-1|^{2}\right)}\ dg&=\int_{V}e^{-t^{2}\left(|e^{X}-1|^{2}+|e^{-X}-1|^{2}\right)}\left|\det\left(\tfrac{1-e^{-\operatorname{ad}(X)}}{\operatorname{ad}(X)}\right)\right| \ dX
\end{align*}
There exists a constant $C>0$ with 
\[ \left|\det\left(\tfrac{1-e^{-\operatorname{ad}(X)}}{\operatorname{ad}(X)}\right)\right|\geq C, \ \ \ \forall X\in V.\]
Hence 
\begin{align*}
\int_{G}e^{-t^{2}\left(|g-1|^{2}+|g^{-1}-1|^{2}\right)}\ dg\geq C \int_{V}e^{-t^{2}\left(|e^{X}-1|^{2}+|e^{-X}-1|^{2}\right)}\ dX
\end{align*}
There exists a constant $C^{\prime}>0$ such that 
\[|e^{X}-1|^{2}+|e^{-X}-1|^{2}\leq C^{\prime} |X|^{2}, \ \ \ \forall X\in V.
\]
Thus
\begin{align*}
\int_{V}e^{-t^{2}\left(|e^{X}-1|^{2}+|e^{-X}-1|^{2}\right)}\ dX &\geq \int_{V}e^{-t^{2}C^{\prime}|X|^{2}}\ dX \\
&= \int_{V}e^{-C^{\prime}|tX|^{2}}\ dX \\
&= \frac{1}{t^{\operatorname{dim}\frak{h}}}\int_{V}e^{-C^{\prime}|X|^{2}}\ dX 
\end{align*}
Therefore
\[\int_{H}e^{-t^{2}\left(|g-1|^{2}+|g^{-1}-1|^{2}\right)}\ dg\geq \int_{V}e^{-C^{\prime}|X|^{2}}\ dX \geq
C_{1}t^{-\operatorname{dim}\frak{h}}
\]
with $C_{1}=C\int_{V}e^{-C^{\prime}|X|^{2}}\ dX  <\infty$. \\
Let $U$ be a neighborhood of $1$ in $G$, there exists a constant $R>0$ such that
\[|g-1|^{2}+|g^{-1}-1|^{2}\geq R, \ \ \ \forall g\in G\backslash U.
\]
Hence 
\[e^{-\frac{1}{2}t^{2}\left(|g-1|^{2}+|g^{-1}-1|^{2}\right)}\leq e^{-\frac{1}{2}t^{2}R}, \ \ \ \forall g\in G\backslash U. \]
Therefore
\begin{align*}
&\int_{G\backslash U}e^{-t^{2}\left(|g-1|^{2}+|g^{-1}-1|^{2}\right)}dg \\
&=\int_{G\backslash U}e^{-\frac{1}{2}t^{2}\left(|g-1|^{2}+|g^{-1}-1|^{2}\right)}
e^{-\frac{1}{2}t^{2}\left(|g-1|^{2}+|g^{-1}-1|^{2}\right)}dg \\
& \leq e^{-\frac{1}{2}t^{2}R}\int_{G\backslash
U}e^{-\frac{1}{2}t^{2}\left(|g-1|^{2}+|g^{-1}-1|^{2}\right)}dg \\
&\leq e^{-\frac{1}{2}t^{2}R}\int_{G\backslash
U}e^{-\frac{1}{2}\left(|g-1|^{2}+|g^{-1}-1|^{2}\right)}dg \\
&= C_{2} e^{-\frac{1}{2}t^{2}R}
\end{align*}
with $C_{2}=\int_{G\backslash
U}e^{-\frac{1}{2}\left(|g-1|^{2}+|g^{-1}-1|^{2}\right)}dg <\infty$. Hence
\begin{align*}
\int_{G\backslash U} \varphi_{t}(g)dh &\leq e^{-\frac{1}{2}t^{2}R}t^{-\operatorname{dim}\frak{g}}\frac{C_{2}}{C_{1}}.
\end{align*}
The expression on the right hand side tends to $0$ as $t$ tends to infinity.
\end{proof}

\begin{lemma}\label{av}
\[\Pi(\mathcal{A}(G))E\subset E^{\omega} \] 
\end{lemma}
\begin{proof}
Let $f\in \mathcal{A}\left(G\right),v\in E$. Then the orbit map $\gamma_{\Pi(f)v}$ is given by 
\begin{align*}
\gamma_{\Pi(f)v}(g)& = \pi(g)\int_{H}f(x)\pi(x)v \ d\mu(x) \\
&= \int_{H}f(x)\pi(gx)v \ d\mu(x) \\
&= \int_{H}f(g^{-1}x)\pi(x)v \ d\mu(x) \\
&=\pi\bigl(L_{g}(f)\bigr)v.
\end{align*}
Hence the orbit map is equal to to the composition
\begin{align*}
G\rightarrow \mathcal{R}(G) \rightarrow  E 
\end{align*}
Here the first arrow denotes the map $g\mapsto L_{g}(f)$ and the second
the map $\varphi\mapsto \Pi(\varphi)v$.
The first map in this composition is analytic and the last is linear. Hence
the whole map is an analytic map from $G$ to $E$.
\end{proof}
\begin{thm}
For every real linear algebraic group $G$ exists an analytic 
Dirac sequence, i.e a Dirac sequence which members are elements of $\mathcal{A}(G)$.
\end{thm}
\begin{proof}
The sequence of functions $(\varphi_{t})_{t\geq1}$ on $G$ provides 
a Dirac sequence, as we have seen in Proposition \ref{d}. 
\end{proof}
\begin{cor}Let $(\pi,E)$ be a $F$-representation of a real linear algebraic group $G$ on a Fr\'{e}chet space $E$. Then
the space $E^{\omega}$ of analytic vectors is dense in $E$.
\end{cor}
\begin{proof}
Let $v\in E$ and let $(\varphi_{t})_{t\geq1}$ be an analytic Dirac sequence.
Then $\pi(\varphi_{t})v$ is, according to Lemma \ref{av}, a sequence of analytic vectors 
which tends to $v$ in $E$.
\end{proof}
\begin{ack}
I want to express my gratitude to my advisor Bernhard Kr\"{o}tz for his generous help.
\end{ack}

\
%\section{Appendix}
%\subsection{Complexification of a linear Lie group}
%Let $H$ be a closed subgroup of $\mathrm{GL}_{n}(\mathbb{R})$ with Lie algebra 
%$\mathfrak{h}$. We denote by $\mathfrak{h}_{\mathbb{C}}$ its complexification.
%There exists an analytic subgroup $H_{\mathbb{C}}^{\circ}$ of $\mathrm{GL}_{n}(\mathbb{C})$
%with Lie algebra $\mathfrak{h}_{\mathbb{C}}$. Then the complexification $H_{\mathbb{C}}$ of $H$ is the 
%group $HH_{\mathbb{C}}^{\circ}$. 

\end{document}